\documentclass[12pt,a4paper]{article}

\usepackage[utf8]{inputenc}
\usepackage[english]{babel}
\usepackage{geometry} 
\usepackage{titlesec}
\usepackage{hyperref}
\usepackage{yfonts}
\usepackage{amssymb}
\usepackage{amsmath}
\usepackage{amsthm}
\usepackage{comment}
\usepackage{xcolor}

\newcommand{\metrica}[1]{g\left(#1\right)}

\newcommand{\ric}[0]{\textup{ric}}
\newcommand{\ii}{\textbf{i}}
\renewcommand{\j}{\textbf{j}}
\renewcommand{\k}{\textbf{k}}
\newcommand{\im}{\textup{Im}}
\newcommand{\A}{$\mathcal{A}$ }

\title{Inaudibility of naturally reductive property}
\author{Teresa Arias-Marco\footnote{ORCID: 0000-0003-0984-0367;\ email: ariasmarco@unex.es} and Jos\'e-Manuel Fern\'andez-Barroso\footnote{ORCID: 0000-0003-3864-9967;\ email: ferbar@unex.es}}
\date{Universidad de Extremadura, Departamento de Matemáticas, Badajoz, Spain.}

\begin{document}
\newtheorem{theorem}{Theorem}[section]
\newtheorem{proposition}[theorem]{Proposition}
\newtheorem{corollary}[theorem]{Corollary}
\newtheorem{lemma}[theorem]{Lemma}

\newtheorem*{main}{Main Theorem}
\newtheorem*{maincor1}{Main corollary 1}
\newtheorem*{maincor2}{Main corollary 2}

\theoremstyle{definition}
\newtheorem{definition}[theorem]{Definition}

\theoremstyle{definition}
\newtheorem{remark}[theorem]{Remark}

\theoremstyle{definition}
\newtheorem{example}[theorem]{Example}

\theoremstyle{definition}
\newtheorem*{notation}{Notation}

\maketitle

\begin{abstract}
In this paper, we use a characterization of naturally reductive 2-step nilponent Lie groups via Ambrose-Singer's homogeneous structures to prove that one cannot determine if a closed Riemannian manifold is naturally reductive using the information encoded in the spectrum of the Laplace-Beltrami operator. To do that, we consider a new isospectral pair of 2-step nilmanifolds of dimension 9 such that one of them is naturally reductive and the other is not.
\end{abstract}



\section*{Introduction and preliminaries}

Let $(M,g)$ denotes a Riemannian manifold and $R(X,Y)=[\nabla_X,\nabla_Y]-\nabla_{[X,Y]}$, $X,Y\in\mathfrak{X}(M)$, its Riemannian curvature operator, where $\nabla$ denotes the Levi-Civita connection. The Ricci operator $\rho$ of $(M,g)$ is the trace in $Y$ of the Jacobi operator $R_Y(X):=R(X,Y)Y$. A Riemannian manifold is \textit{homogeneous} if its group of isometries, $I(M)$, acts transitively on $M$. In particular, following Wilson \cite{W.82}, if $I(M)$ contains a nilpotent Lie subgroup acting transitively on $M$, then $M$ is called a \textit{homogeneous nilmanifold}.  It is well known that locally symmetric spaces are characterize by $\nabla R=0$. Moreover, a \textit{symmetric-like} property of a Riemannian manifold $(M,g)$ is a geometric property that every locally symmetric space satisfies. It is still open if the locally symmetry is audible or not. We will approach to this question studying the audibility of the naturally reductive property. A Riemannian manifold is said naturally reductive if every geodesic in $M$ is the orbit of a one-parameter subgroup of $I(M)$ generated by some vector in the tangent space of $M$ at the origin. Then, every Riemannian manifold satisfying the naturally reductive property is homogeneous. Gordon proved in \cite{G.85} that every naturally reductive nilmanifold is at most 2-step nilpotent. Moreover, naturally reductive Riemannian manifolds are classified up to dimension eight (see \cite{TV.83,KV.83,KV.85,AFF.15,S.20}).

On the other hand, Gray introduced a symmetric-like property, namely Type \A in \cite{G.78}. A Riemannian manifold $(M,g)$ satisfies that property if its Ricci tensor, $\textup{ric}(X,Y):=\metrica{\rho(X),Y}$, is cyclic parallel, this is
\begin{equation}\label{eq:typeA}
    (\nabla_X\ric)(X,X)=0
\end{equation}
for $X,Y\in\mathfrak{X}(M)$. Kowalski and Vanhecke proved in \cite{KV.84} that every Riemannian manifold satisfying the naturally reductive property, also satisfies the D'Atri property. Thus, in particular it also has cyclic parallel Ricci tensor.

In Section \ref{sec:geometria2step} we focus on characterize the naturally reductive property on some kind of 2-step nilpotent Lie groups. Our result is an adaptation of the one given by Gordon in \cite{G.85} and Lauret in \cite{L.98.2}. To obtain it, we first need to characterize the Type \A property on 2-step nilpotent Lie groups. Then, we use it together with the characterization of the naturally reductive property given by Tricerri and Vanhecke in \cite{TV.83} via homogeneous structures. Ambrose and Singer \cite{AS.58} introduced homogeneous structures in order to characterize homogeneous Riemannian manifolds using a local condition which must be satisfied at all points. Tricerri and Vanhecke \cite{TV.83} studied the different classes of homogeneity classifying the homogeneous structures. Recently, Calviño-Louzao, Ferreiro-Subrido, Garc\'ia-R\'io and V\'azquez-Lorenzo \cite{CFGV.23} have determined all the possible homogeneous structures on non-symmetric three-dimensional Riemannian Lie groups. Moreover, Castrill\'on and Calvaruso \cite{CC.19} have extended the Ambrose-Singer result to pseudo-Riemannian manifolds. More precisely, we say that a $(2,1)$-tensor, $T$, is a \textit{homogeneous structure} on a Riemannian manifold $(M,g)$ if
\begin{equation}\label{theo:AS-estructura(1)}
(\nabla_V R)(X,Y)Z=T_V(R(X,Y)Z)-R(T_VX,Y)Z-R(X,T_VY)Z-R(X,Y)(T_VZ),
\end{equation}
\begin{equation}\label{theo:AS-estructura(2)}
(\nabla_XT)_Y=[T_X,T_Y]-T_{T_XY},
\end{equation}
and
\begin{equation}\label{theo:AS-estructura(3)}
g(T_XY,Z)+g(Y,T_XZ)=0
\end{equation}
for $X,Y,Z\in\mathfrak{X}(M)$. Note that, if $\tilde{\nabla}=\nabla-T$, then \eqref{theo:AS-estructura(1)}, \eqref{theo:AS-estructura(2)} and \eqref{theo:AS-estructura(3)} are equivalent to
$$
\tilde{\nabla}R=0,\quad\quad\tilde{\nabla}T=0,\quad\quad\tilde{\nabla}g=0,
$$
respectively. This is, $R$, $T$ and $g$ are parallel with respect to $\tilde{\nabla}$.

\begin{theorem}
A Riemannian manifold $(M,g)$ is locally homogeneous if and only if it has a homogeneous structure, $T$. In particular, $(M,g)$ is naturally reductive if and only if $T$ additionally satisfies
\begin{equation}\label{pro:NR}
g(T_XY,Z)+g(T_YX,Z)=0
\end{equation}
for $X,Y,Z\in\mathfrak{X}(M)$, or equivalently $T_XX=0$.
\end{theorem}

Now, we obtain a short proof of the fact that every naturally reductive Riemannian manifold is a Type \A manifold. The equation obtained in this proof will be useful to obtain Theorem \ref{theo:caracterizacionNR}. 

\begin{proposition}
Every naturally reductive Riemannian manifold is of Type $\mathcal{A}$.
\end{proposition}
\begin{proof}
Considering on \eqref{theo:AS-estructura(1)} the trace of the curvature tensor in order to obtain the Ricci tensor, every locally homogeneous space also satisfies
\begin{equation}\label{eq:contraccionAS}
(\nabla_X\ric)(Y,Z)=-\ric(T_XY,Z)-\ric(Y,T_XZ)
\end{equation}
for $X,Y,Z\in\mathfrak{X}(M)$. Now, suppose that $(M,g)$ is naturally reductive, then $T_XX=0$ for $X\in\mathfrak{X}(M)$, and the above equation gives \eqref{eq:typeA}.
\end{proof}

Two closed Riemannian manifolds, $(M,g)$ and $(M',g')$, are said \textit{isospectral} for the Laplace-Beltrami operator, $\Delta$ and $\Delta'$ respectively, if the eigenvalues of this operator acting on functions of $M$ and $M'$ respectively, coincide, counting multiplicities, that is
    $$
    \text{Spec}(M,\Delta)=\text{Spec}(M',\Delta').
    $$
Given two isospectral Riemannian manifolds, their geometric properties could differ between them. That properties are named \textit{inaudible}. Otherwise, if a geometric property could be recovered using only the information given by the spectrum of the Laplace-Beltrami operator, it is said to be \textit{audible}. Thus, in order to prove the inaudibility of a property, it is enough to find a pair of isospectral manifolds such that one of them satisfies that property and the other does not.

Since Kac opened the audibility problem in \cite{K.66}, there have been proved to be audible or inaudible different geometric properties using different eigenvalue problems and isospectral Riemannian manifolds. In particular, Gordon, Gornet, Schueth, Webb and Wilson constructed an isospectral pair of closed Riemannian manifolds in \cite{GGSWW.98} having different scalar curvatures. Independently, Szab\'o in \cite{Sz.99} used the same isospectral pair to prove that the local homogeneity cannot be determined from the eigenvalues of the Laplace-Beltrami operator. Using Szab\'o's examples, first author and Schueth proved in \cite{AS.10} the inaudibility of the weakly locally symmetry, the D'Atri property, the probabilistic commutativity, the $\mathfrak{C}$ property and the type $\mathcal{A}$ property. Moreover, we proved in \cite{AF.19} the inaudibility of the $k$-D'Atri property for all values of $k$, and in particular for each value of $k$, using the same closed isospectral Riemannian manifolds. Moreover, using different pairs of isospectral closed Riemannian manifolds, first author and Schueth proved in \cite{AS.17} the inaudibility of some sixth order curvature invariants.

In Section \ref{sec:ejemplos}, we present a new pair of isospectral 2-step nilmanifolds. Then, in Section \ref{sec:inaudibilidadNR}, we use this pair to prove the inaudibility of the naturally reductive property on compact Riemannian manifolds.

\section{The geometry of 2-step nilpotent Lie groups}\label{sec:geometria2step}

This section is dedicated to explore the geometry of 2-step nilpotent Lie groups. It was well studied by Eberlein in \cite{E.94}, but we include it for the sake of completeness. Moreover, we obtain some new results involving the covariant derivative of the Ricci tensor in order to characterize the 2-step nilpotent Lie groups of Type \A. Finally, we use it to characterize some kind of 2-step nilpotent Lie groups which are naturally reductive.

\subsection{Preliminaries of 2-step nilpotent Lie groups}

Let $\mathfrak{v}$ and $\mathfrak{z}$ real vector spaces perpendicular with respect to an inner product, $g$, and denote $\mathfrak{v}\oplus\mathfrak{z}$ as $\mathfrak{n}$ endowed with this inner product. From now on, we denote every vector field in the Lie algebra by $X=X^\mathfrak{v}+X^\mathfrak{z}\in\mathfrak{n}$ where $X^\mathfrak{v}\in\mathfrak{v}$ and $X^\mathfrak{z}\in\mathfrak{z}$. Let a linear map $j:\mathfrak{z}\to\mathfrak{so(v)}$, and define a Lie bracket on $\mathfrak{n}$ by
\begin{equation}\label{eq:relacion-j-corchete}
\metrica{[X^\mathfrak{v},Y^\mathfrak{v}]^j,Z^\mathfrak{z}}=\metrica{j_{Z^\mathfrak{z}}X^\mathfrak{v},Y^\mathfrak{v}}
\end{equation}
for $X^\mathfrak{v},Y^\mathfrak{v}\in\mathfrak{v}$ and $Z^\mathfrak{z}\in\mathfrak{z}$.

This Lie bracket $[\cdot,\cdot]^j$ makes $\mathfrak{n}$ a 2-step nilpotent Lie algebra, this is $[\mathfrak{n},\mathfrak{n}]^j\subseteq\mathfrak{z}$ and $[\mathfrak{n},\mathfrak{z}]^j=0$. We write $[\cdot,\cdot]$ instead of $[\cdot,\cdot]^j$ if there is no risk of confusion with the $j$ map associated to the bracket. Let's denote by $\mathfrak{n}(j)$ the pair $(\mathfrak{n},j)$ formed by the 2-step nilpotent Lie algebra $\mathfrak{n}$ and the fixed $j$ map, and by $(N(j),g)$ the 2-step nilpotent Lie group whose Lie algebra is $\mathfrak{n}(j)$ with the left-invariant Riemannian metric induced by $g$, which we also name by $g$. Note that the exponential map $\exp^j:\mathfrak{n}(j)\to N(j)$ is a diffeomorphism due to $N(j)$ is simply connected and nilpotent. In order to simplify the notation, sometimes we only write $\exp$ instead of $\exp^j$.

As a Riemannian manifold, the Levi-Civita connection and the Riemannian curvature operator of $(N(j),g)$ are given, respectively, by
\begin{equation}\label{Eq:levi-civita-2step}
\nabla_VX=-\frac{1}{2}j_{V^\mathfrak{z}}X^\mathfrak{v}-\frac{1}{2}j_{X^\mathfrak{z}}V^\mathfrak{v}+\frac{1}{2}[V^\mathfrak{v},X^\mathfrak{v}],
\end{equation}
and
$$
\begin{aligned}
R(X,Y)Z=&\nabla_X\nabla_YZ-\nabla_Y\nabla_XZ-\nabla_{[X,Y]}Z\\
=&\frac{1}{2}j_{[X^\mathfrak{v},Y^\mathfrak{v}]}Z^\mathfrak{v}+\frac{1}{4}j_{[X^\mathfrak{v},Z^\mathfrak{v}]}Y^\mathfrak{v}-\frac{1}{4}j_{[Y^\mathfrak{v},Z^\mathfrak{v}]}X^\mathfrak{v}+\\
&+\frac{1}{4}j_{X^\mathfrak{z}}j_{Z^\mathfrak{z}}Y^\mathfrak{v}-\frac{1}{4}j_{Y^\mathfrak{z}}j_{Z^\mathfrak{z}}X^\mathfrak{v}+\frac{1}{4}j_{X^\mathfrak{z}}j_{Y^\mathfrak{z}}Z^\mathfrak{v}-\frac{1}{4}j_{Y^\mathfrak{z}}j_{X^\mathfrak{z}}Z^\mathfrak{v}-\\
&-\frac{1}{4}\left[X^\mathfrak{v},j_{Y^\mathfrak{z}}Z^\mathfrak{v}\right]+\frac{1}{4}\left[Y^\mathfrak{v},j_{X^\mathfrak{z}}Z^\mathfrak{v}\right]-\frac{1}{4}\left[j_{Z^\mathfrak{z}}X^\mathfrak{v},Y^\mathfrak{v}\right]-\frac{1}{4}\left[X^\mathfrak{v},j_{Z^\mathfrak{z}}Y^\mathfrak{v}\right].
\end{aligned}
$$
Then, the Jacobi operator $R_V(X):=R(X,V)V$ is
$$
\begin{aligned}
R_V(X)=&\frac{3}{4}j_{[X^\mathfrak{v},V^\mathfrak{v}]}V^\mathfrak{v}+\frac{1}{2}j_{X^\mathfrak{z}}j_{V^\mathfrak{z}}V^\mathfrak{v}-\frac{1}{4}j_{V^\mathfrak{z}}j_{X^\mathfrak{z}}V^\mathfrak{v}-\frac{1}{4}j_{V^\mathfrak{z}}j_{V^\mathfrak{z}}X^\mathfrak{v}-\\
&-\frac{1}{2}[X^\mathfrak{v},j_{V^\mathfrak{z}}V^\mathfrak{v}]+\frac{1}{4}[V^\mathfrak{v},j_{X^\mathfrak{z}}V^\mathfrak{v}]-\frac{1}{4}[j_{V^\mathfrak{z}}X^\mathfrak{v},V^\mathfrak{v}].
\end{aligned}
$$

Let $\{v_i\}, i=1,\dots,n$ be an orthonormal basis of $\mathfrak{v}$ and $\{z_k\}$,  $k=1,\dots,m$ be an orthonormal basis of $\mathfrak{z}$, where $n=\dim\mathfrak{v}$ and $m=\dim\mathfrak{z}$. So, the Ricci operator and the Ricci tensor are given, respectively, by
$$
\begin{aligned}
\rho(X)&=\sum_{i=1}^nR_{v_i}X+\sum_{k=1}^mR_{z_k}X\\
&=\frac{1}{2}\sum_{k=1}^mj_{z_k}j_{z_k}X^\mathfrak{v}+\frac{1}{4}\sum_{i=1}^n\left[v_i,j_{X^\mathfrak{z}}v_i\right],
\end{aligned}
$$
and
$$
\textup{ric}(X,Y)=\metrica{\rho(X),Y}=\frac{1}{2}\sum_{k=1}^m\metrica{j_{z_k}j_{z_k}X^\mathfrak{v},Y^\mathfrak{v}}+\frac{1}{4}\sum_{i=1}^n\metrica{\left[v_i,j_{X^\mathfrak{z}}v_i\right],Y^\mathfrak{z}}.
$$

In order to simplify the notation, we introduce the maps $J:\mathfrak{v}\to\mathfrak{v}$ and $B:\mathfrak{z}\to\mathfrak{z}$, respectively by
$$
J(X^\mathfrak{v}):=\sum_{k=1}^mj_{z_k}j_{z_k}X^\mathfrak{v},
\quad\text{and}\quad
B(X^\mathfrak{z}):=\sum_{i=1}^n\left[v_i,j_{X^\mathfrak{z}}v_i\right].
$$

These maps are $g$-equivariant in the following sense:
$$
\metrica{J(X^\mathfrak{v}),Y^\mathfrak{v}}=\sum_{k=1}^m\metrica{j_{z_k}j_{z_k}X^\mathfrak{v},Y^\mathfrak{v}}=\sum_{k=1}^m\metrica{X^\mathfrak{v},j_{z_k}j_{z_k}Y^\mathfrak{v}}=\metrica{X^\mathfrak{v},J(Y^\mathfrak{v})},
$$
and
$$
\begin{aligned}
\metrica{B(X^\mathfrak{z}),Y^\mathfrak{z}}&=\sum_{i=1}^n\metrica{\left[v_i,j_{X^\mathfrak{z}}v_i\right],Y^\mathfrak{z}}=\sum_{i=1}^n\metrica{j_{Y^\mathfrak{z}}v_i,j_{X^\mathfrak{z}}v_i}=\sum_{i=1}^n\metrica{j_{X^\mathfrak{z}}v_i,j_{Y^\mathfrak{z}}v_i}\\
&=\sum_{i=1}^n\metrica{\left[v_i,j_{Y^\mathfrak{z}}v_i\right],X^\mathfrak{z}}=\metrica{X^\mathfrak{z},B(Y^\mathfrak{z})}.
\end{aligned}
$$

Thus, we can rewrite the Ricci operator and the Ricci tensor as follows
\begin{equation}\label{Eq:RicciOperator-2step}
\rho(X)=\frac{1}{2}J(X^\mathfrak{v})+\frac{1}{4}B(X^\mathfrak{z}),
\end{equation}
\begin{equation}\label{Eq:RicciTensor-2step}
\textup{ric}(X,Y)=\frac{1}{2}\metrica{J(X^\mathfrak{v}),Y^\mathfrak{v}}+\frac{1}{4}\metrica{B(X^\mathfrak{z}),Y^\mathfrak{z}}.
\end{equation}

This $J$ map was used by the first author and Schueth in \cite{AS.17} to prove the inaudibility of some sixth order curvature invariants.

\subsection{2-step nilpotent Lie groups of Type \A}

Here, we are going to characterize the 2-step nilpotent Lie groups of Type \A via the $J$ endomorphism.

\begin{proposition}\label{prop:derivadaricci}
The covariant derivative of the Ricci tensor on a 2-step nilpotent Lie group $(N(j),g)$ is given by
\begin{equation}\label{eq:NablaTensorRicci}
\begin{aligned}
(\nabla_X\textup{ric})(Y,Z)=&\frac{1}{4}\metrica{j_{X^\mathfrak{z}}Y^\mathfrak{v},J(Z^\mathfrak{v})}+\frac{1}{4}\metrica{j_{X^\mathfrak{z}}Z^\mathfrak{v},J(Y^\mathfrak{v})}+\\
&+\frac{1}{4}\metrica{j_{Y^\mathfrak{z}}X^\mathfrak{v},J(Z^\mathfrak{v})}+\frac{1}{4}\metrica{j_{Z^\mathfrak{z}}X^\mathfrak{v},J(Y^\mathfrak{v})}-\\
&-\frac{1}{8}\metrica{\left[X^\mathfrak{v},Z^\mathfrak{v}\right],B(Y^\mathfrak{z})}-\frac{1}{8}\metrica{\left[X^\mathfrak{v},Y^\mathfrak{v}\right],B(Z^\mathfrak{z})}
\end{aligned}
\end{equation}
for $X,Y,Z\in\mathfrak{n}$.
\end{proposition}

\begin{proof}
Using \eqref{Eq:levi-civita-2step} and \eqref{Eq:RicciOperator-2step}, it is easy to get
$$
\begin{aligned}
(\nabla_X\rho)(Y)=&\nabla_X(\rho(Y))-\rho(\nabla_XY)\\
=&\frac{1}{4}J(j_{X^\mathfrak{z}}Y^\mathfrak{v})-\frac{1}{4}j_{X^\mathfrak{z}}J(Y^\mathfrak{v})+\frac{1}{4}J(j_{Y^\mathfrak{z}}X^\mathfrak{v})+\\
&+\frac{1}{4}\left[X^\mathfrak{v},J(Y^\mathfrak{v})\right]-\frac{1}{8}j_{B(Y^\mathfrak{z})}X^\mathfrak{v}-\frac{1}{8}B([X^\mathfrak{v},Y^\mathfrak{v}]).
\end{aligned}
$$
Now, the proposition follows using the $g$-equivariant property of $B$ and $J$ on $(\nabla_X\textup{ric})(Y,Z)=\metrica{(\nabla_X\rho)(Y),Z}$.
\end{proof}

\begin{theorem}\label{theo:caracterizacionA}
A 2-step nilpotent Lie group $(N(j),g)$ is a Type $\mathcal{A}$ manifold if and only if $J\circ j_{X^\mathfrak{z}}$ is a skew-symmetric endomorphism for every $X^\mathfrak{z}\in\mathfrak{z}$.
\end{theorem}
\begin{proof}
By Proposition \ref{prop:derivadaricci}, we get
$$
(\nabla_X\textup{ric})(X,X)=\metrica{J(j_{X^\mathfrak{z}}X^\mathfrak{v}),X^\mathfrak{v}}
$$
for $X\in\mathfrak{n}$.

Then, the Type \A condition \eqref{eq:typeA} is equivalent to
$$
\metrica{J(j_{X^\mathfrak{z}}(Y^\mathfrak{v}+Z^\mathfrak{v})),(Y^\mathfrak{v}+Z^\mathfrak{v})}=0
$$
for each $Y^\mathfrak{v},Z^\mathfrak{v}\in\mathfrak{v}$, $X^\mathfrak{z}\in\mathfrak{z}$. Using the linearity of $j$ and $J$, we get the skew-symmetry of $J\circ j_{X^\mathfrak{z}}$, this is
\begin{equation}\label{eq:skew-Joj}
\metrica{J(j_{X^\mathfrak{z}}Y^\mathfrak{v}),Z^\mathfrak{v}}=-\metrica{Y^\mathfrak{v},J(j_{X^\mathfrak{z}}Z^\mathfrak{v})}.
\end{equation}

The opposite direction is obtained taking $Y^\mathfrak{v}=Z^\mathfrak{v}$ in \eqref{eq:skew-Joj}.
\end{proof}

A particular class of 2-step nilpotent Lie groups which is well studied is the class of generalized Heisenberg groups. They were introduced by Kaplan in \cite{K.81}. A 2-step nilpotent Lie algebra, $\mathfrak{n}(j)$, with inner product $g$, satisfying
$$
  (j_{Z^\mathfrak{z}})^2=-\|Z^\mathfrak{z}\|^2\cdot \textup{Id}_\mathfrak{v}
$$
for every $Z^\mathfrak{z}\in\mathfrak{z}$, is named \textit{generalized Heisenberg algebra}. The corresponding simply connected Lie group, $N(j)$ with the induced left-invariant Riemannian metric, is named \textit{generalized Heisenberg group}, and the $j$ map is called a \textit{map of Heisenberg type}. Note that, when $\dim\mathfrak{z}=1$, the Lie algebra is the classical Heisenberg algebra.

\begin{corollary}\label{coro:typeA}
Every 2-step nilpotent Lie group $(N(j),g)$ such that $J=C\cdot \textup{Id}_\mathfrak{v}$, for some $C\in\mathbb{R}$ and all $X^\mathfrak{v}\in\mathfrak{v}$, is a Type $\mathcal{A}$ manifold.
\end{corollary}
\begin{proof}
This result follows from the fact that $J\circ j_{X^\mathfrak{z}}=C\cdot j_{X^\mathfrak{z}}$ and $j_{X^\mathfrak{z}}$ is a skew-symmetric endomorphism. Thus, $J\circ j_{X^\mathfrak{z}}$ is also a skew-symmetric endomorphism.
\end{proof}

\begin{example}
If $(N(j),g)$ is a generalized Heisenberg group, we get 
$$
J(X^\mathfrak{v})=\sum_{k=1}^mj_{z_k}j_{z_k}X^\mathfrak{v}=\sum_{k=1}^m-\|z_k\|^2\cdot  X^\mathfrak{v}=-\dim\mathfrak{z}\cdot  X^\mathfrak{v}.
$$
This means that every generalized Heisenberg group satisfies $J(X^\mathfrak{v})=C\cdot X^\mathfrak{v}$, with constant $C=-\dim\mathfrak{z}$, and they are Type $\mathcal{A}$ manifolds by Corollary \ref{coro:typeA}. 
Moreover, for each $X^\mathfrak{z}\in\mathfrak{z}$, every generalized Heisenberg group satisfies $B(X^\mathfrak{z})=\dim\mathfrak{v}\cdot X^\mathfrak{z}$, due to
$$
\begin{aligned}
\metrica{B(X^\mathfrak{z}),Y^\mathfrak{z}}&=\sum_{i=1}^n\metrica{\left[v_i,j_{X^\mathfrak{z}}v_i\right],Y^\mathfrak{z}}=\sum_{i=1}^n\metrica{j_{Y^\mathfrak{z}}v_i,j_{X^\mathfrak{z}}v_i}\\
&=\sum_{i=1}^n\|v_i\|^2\metrica{X^\mathfrak{z},Y^\mathfrak{z}}=\dim\mathfrak{v}\cdot\metrica{X^\mathfrak{z},Y^\mathfrak{z}}
\end{aligned}
$$
for every $Y^\mathfrak{z}\in\mathfrak{z}$.
\end{example}

\begin{theorem}
A 2-step nilpotent Lie group $(N(j),g)$ with $J=C\cdot\textup{Id}_\mathfrak{v}$ and $B=D\cdot\textup{Id}_\mathfrak{z}$ has parallel Ricci tensor if and only if $D=2C$.
\end{theorem}
\begin{proof}
If $J=C\cdot\textup{Id}_\mathfrak{v}$ and $B=D\cdot\textup{Id}_\mathfrak{z}$, \eqref{eq:NablaTensorRicci} becomes
\begin{equation}\label{eq:NablaTensorRicci-JBconstant}
(\nabla_X\textup{ric})(Y,Z)=\frac{2C-D}{8}\cdot\left(\metrica{j_{Y^\mathfrak{z}}X^\mathfrak{v},Z^\mathfrak{v}}+\metrica{[X^\mathfrak{v},Y^\mathfrak{v}],Z^\mathfrak{z}}\right)
\end{equation}
for $X,Y,Z\in\mathfrak{n}$. This vanishes if and only if $2C-D=0$.
\end{proof}
\begin{corollary}
None generalized Heisenberg group has parallel Ricci tensor.
\end{corollary}

Finally, we characterize naturally reductive homogeneous structures of any 2-step nilpotent Lie group under some conditions involving the maps $J$ and $B$. This is an adaptation of the theorem given by Gordon in \cite{G.85} and its alternative proof given by Lauret in \cite{L.98.2}.

\begin{theorem}\label{theo:caracterizacionNR}
Let $(N(j),g)$ be a 2-step nilpotent Lie group with $J=C\cdot\textup{Id}_\mathfrak{v}$ and $B=D\cdot\textup{Id}_\mathfrak{z}$, $C,D\in\mathbb{R}$, and $D\neq2C$. Then, its unique naturally reductive homogeneous structure must be
$$
T_XY=-\frac{1}{2}j_{Y^\mathfrak{z}}X^\mathfrak{v}+\frac{1}{2}j_{X^\mathfrak{z}}Y^\mathfrak{v}+\frac{1}{2}[X^\mathfrak{v},Y^\mathfrak{v}]+\tilde{T}_{X^\mathfrak{z}}Y^\mathfrak{z},
$$
where $\tilde{T}:\mathfrak{z}\times\mathfrak{z}\to\mathfrak{z}$ is a Lie bracket in $\mathfrak{z}$ such that $\tilde{T}_{X^\mathfrak{z}}\in\mathfrak{so(z)}$, and
\begin{equation}
\label{eq:condicion2-theorem} j_{X^\mathfrak{z}}j_{Y^\mathfrak{z}}Z^\mathfrak{v}-j_{Y^\mathfrak{z}}j_{X^\mathfrak{z}}Z^\mathfrak{v}=j_{\tilde{T}_{X^\mathfrak{z}}Y^\mathfrak{z}}Z^\mathfrak{v}
\end{equation}
for each $X,Y,Z\in\mathfrak{n}$.
\end{theorem}

\begin{proof}
Every homogeneous structure $T$ in any 2-step nilpotent Lie group could be written as follows
\begin{equation}\label{eq:EstructuraT-2step-0}
\begin{aligned}
T_XY=&(T_XY)^\mathfrak{v}+(T_XY)^\mathfrak{z}\\
=&(T_{X^\mathfrak{v}}Y^\mathfrak{v})^\mathfrak{v}+(T_{X^\mathfrak{v}}Y^\mathfrak{z})^\mathfrak{v}+(T_{X^\mathfrak{z}}Y^\mathfrak{v})^\mathfrak{v}+(T_{X^\mathfrak{z}}Y^\mathfrak{z})^\mathfrak{v}\\
&+(T_{X^\mathfrak{v}}Y^\mathfrak{v})^\mathfrak{z}+(T_{X^\mathfrak{v}}Y^\mathfrak{z})^\mathfrak{z}+(T_{X^\mathfrak{z}}Y^\mathfrak{v})^\mathfrak{z}+(T_{X^\mathfrak{z}}Y^\mathfrak{z})^\mathfrak{z}.
\end{aligned}
\end{equation}

This homogeneous structure should satisfies \eqref{theo:AS-estructura(3)}, thus $T_X\in\mathfrak{so(n)}$ for every $X\in\mathfrak{n}$.

If $J=C\cdot\textup{Id}_\mathfrak{v}$ and $B=D\cdot\textup{Id}_\mathfrak{z}$, then \eqref{Eq:RicciTensor-2step} becomes
\begin{equation}\label{eq:tensorRicci-JBconstant}
\textup{ric}(X,Y)=\frac{C}{2}\metrica{X^\mathfrak{v},Y^\mathfrak{v}}+\frac{D}{4}\metrica{X^\mathfrak{z},Y^\mathfrak{z}}.
\end{equation}

Using \eqref{eq:EstructuraT-2step-0} and \eqref{eq:tensorRicci-JBconstant}, we obtain
\begin{equation}\label{eq:-ric-ric - JBconstant}
\begin{aligned}
-\textup{ric}(T_XY,Z)-\textup{ric}(Y,T_XZ)=&-\frac{C}{2}\metrica{(T_XY)^\mathfrak{v},Z^\mathfrak{v}}-\frac{D}{4}\metrica{(T_XY)^\mathfrak{z},Z^\mathfrak{z}}\\
&-\frac{C}{2}\metrica{Y^\mathfrak{v},(T_XZ)^\mathfrak{v}}-\frac{D}{4}\metrica{Y^\mathfrak{z},(T_XZ)^\mathfrak{z}}.
\end{aligned}
\end{equation}
Now, in order to satisfy \eqref{eq:contraccionAS}, from \eqref{eq:NablaTensorRicci-JBconstant} and \eqref{eq:-ric-ric - JBconstant}, we obtain that the only non trivial relations are
\begin{equation}\label{eq:a-estructura}
    \left(T_{X^\mathfrak{z}}Y^\mathfrak{z}\right)^\mathfrak{v}=0\text{ for }Z\in\mathfrak{v}, X, Y\in\mathfrak{z},
\end{equation}
\begin{equation}\label{eq:b-estructura}
    \left(T_{X^\mathfrak{z}}Y^\mathfrak{v}\right)^\mathfrak{z}=0\text{ for }Y\in\mathfrak{v}, X, Z\in\mathfrak{z},
\end{equation}
\begin{equation}\label{eq:c-estructura}
    \left(T_{X^\mathfrak{v}}Y^\mathfrak{z}\right)^\mathfrak{v}=-\frac{1}{2}j_{Y^\mathfrak{z}}X^\mathfrak{v}\text{ for }X, Z\in\mathfrak{v}, Y\in\mathfrak{z},
\end{equation}
\begin{equation}\label{eq:d-estructura}
    \left(T_{X^\mathfrak{v}}Y^\mathfrak{v}\right)^\mathfrak{z}=\frac{1}{2}[X^\mathfrak{v},Y^\mathfrak{v}]\text{ for }X, Y\in\mathfrak{v}, Z\in\mathfrak{z}.
\end{equation}
In addition, due to $T$ is a naturally reductive homogeneous structure, from \eqref{pro:NR}, \eqref{eq:b-estructura} and \eqref{eq:c-estructura}, we get
\begin{equation}\label{eq:e-estructura-NR}
    \left(T_{X^\mathfrak{z}}Y^\mathfrak{v}\right)^\mathfrak{v}=\left(-T_{Y^\mathfrak{v}}X^\mathfrak{z}\right)^\mathfrak{v}=\frac{1}{2}j_{X^\mathfrak{z}}Y^\mathfrak{v},
\end{equation}
\begin{equation}\label{eq:f-estructura-NR}
   \left(T_{X^\mathfrak{v}}Y^\mathfrak{z}\right)^\mathfrak{z}=\left(-T_{Y^\mathfrak{z}}X^\mathfrak{v}\right)^\mathfrak{z}=0.
\end{equation}
Thus, \eqref{eq:EstructuraT-2step-0} is simplified using equations \eqref{eq:a-estructura}--\eqref{eq:f-estructura-NR}  as follows
\begin{equation}\label{eq:EstructuraT-2step-1}
T_XY=\bar{T}_{X^\mathfrak{v}}Y^\mathfrak{v}-\frac{1}{2}j_{Y^\mathfrak{z}}X^\mathfrak{v}+\frac{1}{2}j_{X^\mathfrak{z}}Y^\mathfrak{v}+\frac{1}{2}[X^\mathfrak{v},Y^\mathfrak{v}]+\tilde{T}_{X^\mathfrak{z}}Y^\mathfrak{z},
\end{equation}
where $\bar{T}:\mathfrak{v}\times\mathfrak{v}\to\mathfrak{v}$ and $\tilde{T}:\mathfrak{z}\times\mathfrak{z}\to\mathfrak{z}$ are given by $\bar{T}_{X^\mathfrak{v}}Y^\mathfrak{v}=(T_{X^\mathfrak{v}}Y^\mathfrak{v})^\mathfrak{v}$ and $\tilde{T}_{X^\mathfrak{z}}Y^\mathfrak{z}=(T_{X^\mathfrak{z}}Y^\mathfrak{z})^\mathfrak{z}$. Moreover, due to $T$ is naturally reductive, by \eqref{pro:NR} we also know
\begin{equation}\label{eq:condicionEstructuraNR}
    \bar{T}_{X^\mathfrak{v}}Y^\mathfrak{v}=-\bar{T}_{Y^\mathfrak{v}}X^\mathfrak{v}\text{ and }\tilde{T}_{X^\mathfrak{z}}Y^\mathfrak{z}=-\tilde{T}_{Y^\mathfrak{z}}X^\mathfrak{z}
\end{equation}
for $X,Y\in\mathfrak{n}$.

As homogeneous structure, $T$ also satisfy $\tilde{\nabla}g=0$ and $\tilde{\nabla}T=0$, where $\tilde{\nabla}=\nabla-T$. Using \eqref{Eq:levi-civita-2step} and \eqref{eq:EstructuraT-2step-1}, $\tilde{\nabla}$ is given by
\begin{equation*}
    \tilde{\nabla}_XY=-\bar{T}_{X^\mathfrak{v}}Y^\mathfrak{v}-j_{X^\mathfrak{z}}Y^\mathfrak{v}-\tilde{T}_{X^\mathfrak{z}}Y^\mathfrak{z}
\end{equation*}
for $X,Y\in\mathfrak{n}$. Thus,
\begin{equation}\label{eq:nablatildeg}
\begin{aligned}
0=&(\tilde{\nabla}_Xg)(Y,Z)=-\metrica{\tilde{\nabla}_XY,Z}-\metrica{Y,\tilde{\nabla}_XZ}\\
=&\metrica{\bar{T}_{X^\mathfrak{v}}Y^\mathfrak{v},Z^\mathfrak{v}}+\metrica{\tilde{T}_{X^\mathfrak{z}}Y^\mathfrak{z},Z^\mathfrak{z}}+\metrica{Y^\mathfrak{v},\bar{T}_{X^\mathfrak{v}}Z^\mathfrak{v}}+\metrica{Y^\mathfrak{z},\tilde{T}_{X^\mathfrak{z}}Z^\mathfrak{z}},
\end{aligned}
\end{equation}
and
\begin{equation}\label{eq:nablatildeT}
\begin{aligned}
0=&(\tilde{\nabla}_XT)_YZ=\tilde{\nabla}_X(T_YZ)-T_{\tilde{\nabla}_XY}Z-T_Y(\tilde{\nabla}_XZ)\\
=&-\bar{T}_{X^\mathfrak{v}}\bar{T}_{Y^\mathfrak{v}}Z^\mathfrak{v}+\frac{1}{2}\bar{T}_{X^\mathfrak{v}}j_{Z^\mathfrak{z}}Y^\mathfrak{v}-\frac{1}{2}\bar{T}_{X^\mathfrak{v}}j_{Y^\mathfrak{z}}Z^\mathfrak{v}-j_{X^\mathfrak{z}}\bar{T}_{Y^\mathfrak{v}}Z^\mathfrak{v}\\
&+\frac{1}{2}j_{X^\mathfrak{z}}j_{Z^\mathfrak{z}}Y^\mathfrak{v}-\frac{1}{2}j_{X^\mathfrak{z}}j_{Y^\mathfrak{z}}Z^\mathfrak{v}-\frac{1}{2}\tilde{T}_{X^\mathfrak{z}}([Y^\mathfrak{v},Z^\mathfrak{v}])-\tilde{T}_{X^\mathfrak{z}}\tilde{T}_{Y^\mathfrak{z}}Z^\mathfrak{z}\\
&+\bar{T}_{\bar{T}_{X^\mathfrak{v}}Y^\mathfrak{v}}Z^\mathfrak{v}+\bar{T}_{j_{X^\mathfrak{z}}Y^\mathfrak{v}}Z^\mathfrak{v}-\frac{1}{2}j_{Z^\mathfrak{z}}\bar{T}_{X^\mathfrak{v}}Y^\mathfrak{v}-\frac{1}{2}j_{Z^\mathfrak{z}}j_{X^\mathfrak{z}}Y^\mathfrak{v}\\
&+\frac{1}{2}j_{\tilde{T}_{X^\mathfrak{z}}Y^\mathfrak{z}}Z^\mathfrak{v}+\frac{1}{2}\left[\bar{T}_{X^\mathfrak{v}}Y^\mathfrak{v},Z^\mathfrak{v}\right]+\frac{1}{2}\left[j_{X^\mathfrak{z}}Y^\mathfrak{v},Z^\mathfrak{v}\right]+\tilde{T}_{\tilde{T}_{X^\mathfrak{z}}Y^\mathfrak{z}}Z^\mathfrak{z}\\
&+\bar{T}_{Y^\mathfrak{v}}\bar{T}_{X^\mathfrak{v}}Z^\mathfrak{v}+\bar{T}_{Y^\mathfrak{v}}j_{X^\mathfrak{z}}Z^\mathfrak{v}-\frac{1}{2}j_{\tilde{T}_{X^\mathfrak{z}}Z^\mathfrak{z}}Y^\mathfrak{v}+\frac{1}{2}j_{Y^\mathfrak{z}}\bar{T}_{X^\mathfrak{v}}Z^\mathfrak{v}\\
&+\frac{1}{2}j_{Y^\mathfrak{z}}j_{X^\mathfrak{z}}Z^\mathfrak{v}+\frac{1}{2}\left[Y^\mathfrak{v},\bar{T}_{X^\mathfrak{v}}Z^\mathfrak{v}\right]+\frac{1}{2}\left[Y^\mathfrak{v},j_{X^\mathfrak{z}}Z^\mathfrak{v}\right]+\tilde{T}_{Y^\mathfrak{z}}\tilde{T}_{X^\mathfrak{z}}Z^\mathfrak{z}
\end{aligned}
\end{equation}
for $X,Y,Z\in\mathfrak{n}$. Therefore, from \eqref{eq:nablatildeg}, the only non trivial relations are:
\begin{itemize}
    \item If $X, Y, Z\in\mathfrak{z}$,
    $$    \metrica{\tilde{T}_{X^\mathfrak{z}}Y^\mathfrak{z},Z^\mathfrak{z}}+\metrica{Y^\mathfrak{z},\tilde{T}_{X^\mathfrak{z}}Z^\mathfrak{z}}=0.
    $$
    \item If $X, Y, Z\in\mathfrak{v}$,
    $$    \metrica{\bar{T}_{X^\mathfrak{v}}Y^\mathfrak{v},Z^\mathfrak{v}}+\metrica{Y^\mathfrak{v},\bar{T}_{X^\mathfrak{v}}Z^\mathfrak{v}}=0.
    $$
\end{itemize}
This means that $\tilde{T}_{X^\mathfrak{z}}\in\mathfrak{so(z)}$ and $\bar{T}_{X^\mathfrak{v}}\in\mathfrak{so(v)}$ for all $X\in\mathfrak{n}$.

From \eqref{eq:nablatildeT}, we obtain the following non trivial and non equivalent relations:
\begin{itemize}
    \item If $X,Y,Z\in\mathfrak{z}$,
    $$
    \begin{aligned}
        0=&(\tilde{\nabla}_{X^\mathfrak{z}}T)_{Y^\mathfrak{z}}Z^\mathfrak{z}=-\tilde{T}_{X^\mathfrak{z}}\tilde{T}_{Y^\mathfrak{z}}Z^\mathfrak{z}+\tilde{T}_{\tilde{T}_{X^\mathfrak{z}}Y^\mathfrak{z}}Z^\mathfrak{z}+\tilde{T}_{Y^\mathfrak{z}}\tilde{T}_{X^\mathfrak{z}}Z^\mathfrak{z}\\
        =&\tilde{T}_{X^\mathfrak{z}}\tilde{T}_{Z^\mathfrak{z}}Y^\mathfrak{z}+\tilde{T}_{Z^\mathfrak{z}}\tilde{T}_{Y^\mathfrak{z}}X^\mathfrak{z}+\tilde{T}_{Y^\mathfrak{z}}\tilde{T}_{X^\mathfrak{z}}Z^\mathfrak{z}.
    \end{aligned}
    $$
    Thus, $\tilde{T}$ is a Lie bracket on $\mathfrak{z}$ due to this relation together with \eqref{eq:condicionEstructuraNR}.
    \item If $X,Y\in\mathfrak{z},Z\in\mathfrak{v}$,
    $$
    0=(\tilde{\nabla}_{X^\mathfrak{z}}T)_{Y^\mathfrak{z}}Z^\mathfrak{v}=-\frac{1}{2}j_{X^\mathfrak{z}}j_{Y^\mathfrak{z}}Z^\mathfrak{v}+\frac{1}{2}j_{\tilde{T}_{X^\mathfrak{z}}Y^\mathfrak{z}}Z^\mathfrak{v}+\frac{1}{2}j_{Y^\mathfrak{z}}j_{X^\mathfrak{z}}Z^\mathfrak{v}.
    $$
    Thus
    \begin{equation}\label{eq:NablaTildeT-2}
j_{X^\mathfrak{z}}j_{Y^\mathfrak{z}}Z^\mathfrak{v}-j_{Y^\mathfrak{z}}j_{X^\mathfrak{z}}Z^\mathfrak{v}=j_{\tilde{T}_{X^\mathfrak{z}}Y^\mathfrak{z}}Z^\mathfrak{v}.
\end{equation}
This relation is equivalent to $\left[j_{X^\mathfrak{z}}Y^\mathfrak{v},Z^\mathfrak{v}\right]+\left[Y^\mathfrak{v},j_{X^\mathfrak{z}}Z^\mathfrak{v}\right]=\tilde{T}_{X^\mathfrak{z}}([Y^\mathfrak{v},Z^\mathfrak{v}])$, due to $\tilde{T}_{X^\mathfrak{z}}\in\mathfrak{so(z)}$ and the properties of $j$.
    \item If $X,Z\in\mathfrak{v},Y\in\mathfrak{z}$,
    \begin{equation}\label{eq:NablaTildeT-3}
    0=-2(\tilde{\nabla}_{X^\mathfrak{v}}T)_{Y^\mathfrak{z}}Z^\mathfrak{v}=\bar{T}_{X^\mathfrak{v}}j_{Y^\mathfrak{z}}Z^\mathfrak{v}-j_{Y^\mathfrak{z}}\bar{T}_{X^\mathfrak{v}}Z^\mathfrak{v}.
    \end{equation}
This relation is equivalent to $\left[\bar{T}_{X^\mathfrak{v}}Y^\mathfrak{v},Z^\mathfrak{v}\right]+\left[Y^\mathfrak{v},\bar{T}_{X^\mathfrak{v}}Z^\mathfrak{v}\right]=0$, due to $\bar{T}_{X^\mathfrak{v}}\in\mathfrak{so(v)}$ and the properties of $j$.
    \item If $X,Y,Z\in\mathfrak{v}$,
    $$
    \begin{aligned}
    0=(\tilde{\nabla}_{X^\mathfrak{v}}T)_{Y^\mathfrak{v}}Z^\mathfrak{v}=&-\bar{T}_{X^\mathfrak{v}}\bar{T}_{Y^\mathfrak{v}}Z^\mathfrak{v}+\bar{T}_{\bar{T}_{X^\mathfrak{v}}Y^\mathfrak{v}}Z^\mathfrak{v}+\bar{T}_{Y^\mathfrak{v}}\bar{T}_{X^\mathfrak{v}}Z^\mathfrak{v}\\
    &+\frac{1}{2}\left[\bar{T}_{X^\mathfrak{v}}Y^\mathfrak{v},Z^\mathfrak{v}\right]+\frac{1}{2}\left[Y^\mathfrak{v},\bar{T}_{X^\mathfrak{v}}Z^\mathfrak{v}\right].
    \end{aligned}
    $$
    Thus, using \eqref{eq:NablaTildeT-3}, this relation is simplified as follows
$$
\bar{T}_{X^\mathfrak{v}}\bar{T}_{Z^\mathfrak{v}}Y^\mathfrak{v}+\bar{T}_{Z^\mathfrak{v}}\bar{T}_{Y^\mathfrak{v}}X^\mathfrak{v}+\bar{T}_{Y^\mathfrak{v}}\bar{T}_{X^\mathfrak{v}}Z^\mathfrak{v}=0.
$$
Which, together with \eqref{eq:condicionEstructuraNR}, mean that $\bar{T}$ is a Lie bracket in $\mathfrak{v}$.
    \item If $Y,Z\in\mathfrak{v},X\in\mathfrak{z}$,
    $$
    \begin{aligned}
    0=(\tilde{\nabla}_{X^\mathfrak{z}}T)_{Y^\mathfrak{v}}Z^\mathfrak{v}=&-j_{X^\mathfrak{z}}\bar{T}_{Y^\mathfrak{v}}Z^\mathfrak{v}-\frac{1}{2}\tilde{T}_{X^\mathfrak{z}}([Y^\mathfrak{v},Z^\mathfrak{v}])+\bar{T}_{j_{X^\mathfrak{z}}Y^\mathfrak{v}}Z^\mathfrak{v}\\
    &+\frac{1}{2}\left[j_{X^\mathfrak{z}}Y^\mathfrak{v},Z^\mathfrak{v}\right]+\bar{T}_{Y^\mathfrak{v}}j_{X^\mathfrak{z}}Z^\mathfrak{v}+\frac{1}{2}\left[Y^\mathfrak{v},j_{X^\mathfrak{z}}Z^\mathfrak{v}\right].
    \end{aligned}
    $$
Using \eqref{eq:NablaTildeT-2} and \eqref{eq:NablaTildeT-3}, this equation is equivalent to
$$
0=\bar{T}_{Z^\mathfrak{v}}j_{X^\mathfrak{z}}Y^\mathfrak{v}=j_{X^\mathfrak{z}}\bar{T}_{Z^\mathfrak{v}}Y^\mathfrak{v}.
$$
Thus, $j_{X^\mathfrak{z}}\bar{T}_{Z^\mathfrak{v}}Y^\mathfrak{v}$ vanishes if and only if, for any non-zero $W^\mathfrak{v}\in\mathfrak{v}$
$$
\metrica{\left[\bar{T}_{Z^\mathfrak{v}}Y^\mathfrak{v},W^\mathfrak{v}\right],X^\mathfrak{v}}=0.
$$
Due to $X^\mathfrak{v}$ is an arbitrary non-zero vector in $\mathfrak{v}$, this equation vanishes if and only if
$$
\bar{T}_{Z^\mathfrak{v}}Y^\mathfrak{v}=0.
$$
This means that $\bar{T}$ is the abelian Lie bracket in $\mathfrak{v}$.
\end{itemize}
\end{proof}

\section{The new isospectral pair}\label{sec:ejemplos}

We present the construction of our pair of isospectral compact 2-step nilmanifolds.  We take $\mathfrak{v}=\mathbb{R}^6$ and $\mathfrak{z}=\mathbb{R}^3$ with their standard inner products, we can identify $\mathfrak{v}$ with $\mathbb{H}^*\oplus\mathbb{H}^*$, and $\mathfrak{z}$ with $\mathbb{H}^*$, where $\mathbb{H}^*$ the pure quaternion numbers spanned by $\{\ii,\j,\k\}$. In order to distinguish each pure quaternion space, let $\{X_{\ii},X_\j,X_\k,Y_{\ii},Y_\j,Y_\k\}$ and $\{Z_{\ii},Z_\j,Z_\k\}$ be orthonormal bases of $\mathfrak{v}$ and $\mathfrak{z}$, respectively. Then, we define the maps $j,j':\mathfrak{z}\to\mathfrak{so(v)}$ such that, for each $C=(c_1,c_2,c_3)=c_1Z_{\ii}+c_2Z_\j+c_3Z_\k\in\mathfrak{z}$ and each $V=(L,R)\in\mathfrak{v}$, are given by 
$$
\begin{aligned}
j_CV&=j_C(L,R)=(\im(C\cdot L),\im(C\cdot R)),\\
j'_CV&=j'_C(L,R)=(\im(R\cdot C),\im(L\cdot C)),
\end{aligned}
$$
where $\im$ denotes the imaginary part of the quaternion and $\cdot$ is the usual multiplication in $\mathbb{H}$. In order to simplify the computations, we are going to use their matrices representations
\begin{equation}\label{def:matriz-j}
j_C= \begin{pmatrix}
0&c_3 &-c_2&0&0&0\\
-c_3&0&c_1&0&0&0\\
c_2&-c_1&0&0&0&0\\
0&0&0&0&c_3 &-c_2\\
0&0&0&-c_3&0&c_1\\
0&0&0&c_2&-c_1&0
\end{pmatrix},
\end{equation}
and
\begin{equation}\label{def:matriz-j'}
j'_C= \begin{pmatrix}
0&0&0&0&-c_3 &c_2\\
0&0&0&c_3&0&-c_1\\
0&0&0&-c_2&c_1&0\\
0&-c_3 &c_2&0&0&0\\
c_3&0&-c_1&0&0&0\\
-c_2&c_1&0&0&0&0
\end{pmatrix}.
\end{equation}

Using the relation \eqref{eq:relacion-j-corchete} between the Lie bracket and the $j$ map on every 2-step nilpotent Lie algebra, we obtain that the only non vanishing Lie brackets defined on $\mathfrak{n}(j)$ and on $\mathfrak{n}(j')$ are given by
\begin{equation}\label{def:corchete-j}
\left\{
\begin{matrix}
    \left[X_{\ii},X_\j\right] &=&Z_\k, && \left[Y_{\ii},Y_\j\right] &=&Z_\k,\\
    \left[X_{\ii},X_\k\right] &=&-Z_\j, && \left[Y_{\ii},Y_\k\right] &=&-Z_\j,\\
    \left[X_\j,X_\k\right] &=&Z_{\ii}, && \left[Y_\j,Y_\k\right] &=&Z_{\ii},
\end{matrix}
\right.
\end{equation}
and
\begin{equation}\label{def:corchete-j'}
\left\{
\begin{matrix}
    \left[X_i,Y_j\right]'  &=&-Z_k, && \left[Y_i,X_j\right]'  &=&Z_k,\\
    \left[X_i,Y_k\right]'  &=&Z_j, && \left[Y_i,X_k\right]'  &=&-Z_j,\\
    \left[X_j,Y_k\right]'  &=&-Z_i, && \left[Y_j,X_k\right]'  &=&Z_i.
\end{matrix}
\right.
\end{equation}

\begin{remark} 
$\mathfrak{n}(j)$ and $\mathfrak{n}(j')$ are not generalized Heisenberg algebras, due to, for each $C=(c_1,c_2,c_3)=c_1Z_{\ii}+c_2Z_\j+c_3Z_\k\in\mathfrak{z}$, one can obtain
\small
$$
\begin{aligned}
(j_{C^\mathfrak{z}} )^2=(j'_{C^\mathfrak{z}} )^2&=\left(
\begin{array}{ccccccccc}
 -c_2^2-c_3^2 & c_1 c_2 & c_1 c_3 & 0 & 0 & 0 \\
 c_1 c_2 & -c_1^2-c_3^2 & c_2 c_3 & 0 & 0 & 0 \\
 c_1 c_3 & c_2 c_3 & -c_1^2-c_2^2 & 0 & 0 & 0 \\
 0 & 0 & 0 & -c_2^2-c_3^2 & c_1 c_2 & c_1 c_3 \\
 0 & 0 & 0 & c_1 c_2 & -c_1^2-c_3^2 & c_2 c_3 \\
 0 & 0 & 0 & c_1 c_3 & c_2 c_3 & -c_1^2-c_2^2 
\end{array}
\right)\\
&\neq-\|C^\mathfrak{z}\|^2\cdot\textup{Id}_\mathfrak{v}.
\end{aligned}
$$\normalsize
\end{remark}

Lauret introduced in \cite{L.99} a particular class of 2-step nilpotent Lie groups namely \textit{modified Heisenberg groups}, as those whose 2-step nilpotent Lie algebras $\mathfrak{n}(j)$ satisfy
    $$
    (j_{Z^\mathfrak{z}})^2=\lambda(Z^\mathfrak{z})\ \textup{Id}_\mathfrak{v},
    $$
for every nonzero $Z^\mathfrak{z}\in\mathfrak{z}$ and some $\lambda(Z^\mathfrak{z})<0$. It is clear that every generalized Heisenberg group is a modified Heisenberg group considering $\lambda(Z^\mathfrak{z})=-\|Z^\mathfrak{z}\|^2$.

\begin{remark}
The Lie groups associated with the presented 2-step nilpotent Lie algebras $\mathfrak{n}(j)$ and $\mathfrak{n}(j')$ are not modified Heisenberg groups, due to
\begin{equation}\label{eq:j^2-base}
\begin{aligned}
(j_{Z_{\ii}} )^2=(j'_{Z_{\ii}} )^2&=\textup{diag}(0,-1,-1,0,-1,-1)\neq\lambda(Z_{\ii})\ \textup{Id}_\mathfrak{v},\\
(j_{Z_\j} )^2=(j'_{Z_\j} )^2&=\textup{diag}(-1,0,-1,-1,0,-1)\neq\lambda(Z_\j)\ \textup{Id}_\mathfrak{v},\\
(j_{Z_\k} )^2=(j'_{Z_\k} )^2&=\textup{diag}(-1,-1,0,-1,-1,0)\neq\lambda(Z_\k)\ \textup{Id}_\mathfrak{v}.
\end{aligned}
\end{equation}
\end{remark}

Now, let's consider
$$\mathcal{M}=\text{Span}_\mathbb{Z}\{X_{\ii},X_\j,X_\k,Y_{\ii},Y_\j,Y_\k\}\text{ and }\mathcal{L}=\text{Span}_\mathbb{Z}\left\{\tfrac{1}{2}Z_{\ii},\tfrac{1}{2}Z_\j,\tfrac{1}{2}Z_\k\right\}$$
lattices in $\mathfrak{v}$ and $\mathfrak{z}$, respectively. Let $\mathcal{R}=\mathcal{M+L}$ be a lattice in $\mathfrak{n}(j)$ and in $\mathfrak{n}(j')$. Then, the quotient of the 2-step nilpotent Lie group $N(j)$ by a cocompact discrete subgroup $\Gamma(j)=\exp^j(\mathcal{R})$, together with a left-invariant Riemannian metric $g(j)$, is a compact 2-step nilmanifold denoted by $N^j=(N(j)/\Gamma(j),g(j))$. The cannonical projection from $N(j)$ to $N(j)/\Gamma(j)$ is a Riemannian covering. Analogously, we also consider the compact 2-step nilmanifold $N^{j'}=(N(j')/\Gamma(j'),g(j'))$.

\begin{proposition}\label{prop:NjNj'isospectral}
$N^j$ and $N^{j'}$ are isospectral for the Laplace-Beltrami operator.
\end{proposition}
\begin{proof}
Following proposition given in its original form on \cite[\hspace{0cm}3.2, 3.7, 3.8]{GW.97} and adapted on \cite[Proposition 2.4]{S.08}, we need to prove:
\begin{itemize}
    \item[(i)] For each $Z\in\mathfrak{z}$, $j_Z,j'_Z:\mathfrak{v}\to\mathfrak{v}$ have the same eigenvalues in $\mathbb{C}$, counting multiplicities.
    \item[(ii)] $[\mathcal{M},\mathcal{M}]^j$ and $[\mathcal{M},\mathcal{M}]^{j'}$ are contained in $2\mathcal{L}$.
    \item[(iii)] Let $\mathcal{L}^*=\{Z\in\mathfrak{z}|\metrica{Z,\mathcal{L}}\subseteq\mathbb{Z}\}$ be the dual lattice of $\mathcal{L}$ in $\mathfrak{z}$. For each $Z\in\mathcal{L}^*$, the lengths of the elements of the lattices $\ker(j_Z)\cap\mathcal{M}$ and $\ker(j'_Z)\cap\mathcal{M}$, counted with multiplicities, coincide.
\end{itemize}

We get (i) because
$$
\text{Spec}(j_Z)=\left\{0,\pm i\|Z\|\right\}=\text{Spec}(j'_Z ).
$$

By definitions \eqref{def:corchete-j} and \eqref{def:corchete-j'}, we obtain (ii).

Now, we prove (iii). For every $C=(c_1,c_2,c_3)=c_1Z_{\ii}+c_2Z_\j+c_3Z_\k\in\mathfrak{z}$, and particularly for every $Z\in(2\mathbb{Z})^3=\mathcal{L}^*$, the kernels of the maps $j_Z$ and $j'_Z$ can be distinguished as follows:
\begin{itemize}
    \item If $c_1\neq0$,
    $$
    \ker(j_C )=\textup{Span}\left\{\left(1,\frac{c_2}{c_1},\frac{c_3}{c_1},0,0,0\right),\left(0,0,0,1,\frac{c_2}{c_1},\frac{c_3}{c_1}\right)\right\}=\ker(j'_C  ).
    $$
    \item If $c_1=0, c_2\neq0$,
    $$
    \ker(j_C )=\textup{Span}\left\{\left(0,1,\frac{c_3}{c_2},0,0,0\right),\left(0,0,0,0,1,\frac{c_3}{c_2}\right)\right\}=\ker(j'_C  ).
    $$
    \item If $c_1=0, c_2=0, c_3\neq0$,
    $$
    \ker(j_C )=\textup{Span}\left\{\left(0,0,1,0,0,0\right),\left(0,0,0,0,0,1\right)\right\}=\ker(j'_C  ).
    $$
\end{itemize}
In any case, $\ker(j_C)=\ker(j'_C)$, and this means that the lattices $\ker(j_C )\cap\mathcal{M}$ and $\ker(j'_C  )\cap\mathcal{M}$ are the same.
\end{proof}

\begin{proposition}\label{prop:NjNj'noisometricas}
$N^j$ and $N^{j'}$ are not locally isometric.
\end{proposition}
\begin{proof}
Due to the result given by Wilson in \cite{W.82}, it is enough to prove that $\mathfrak{n}(j)$ and $\mathfrak{n}(j')$ are not isomorphic. We can write the following subspaces
$$
\begin{array}{rcl}
    \mathfrak{v}_X&=&\textup{Span}\{X_{\ii},X_\j,X_\k\},\\
    \mathfrak{v}_Y&=&\textup{Span}\{Y_{\ii},Y_\j,Y_\k\},\\
    \mathfrak{v}_{ab}&=&\textup{Span}\{X_a,Y_b\},\ a,b\in\{\ii,\j,\k\}.
\end{array}
$$
Then, $\mathfrak{n}(j')$ has two 6-dimensional abelian subspaces, $\mathfrak{a}'_1=\mathfrak{v}_X\oplus\mathfrak{z}$ and $\mathfrak{a}'_2=\mathfrak{v}_Y\oplus\mathfrak{z}$, while $\mathfrak{n}(j)$ has nine 5-dimensional abelian subspaces, $\mathfrak{a}_{ab}=\mathfrak{v}_{ab}\oplus\mathfrak{z},\ a,b\in\{\ii,\j,\k\}$.
\end{proof}

\section{The geometry of $N^j$ and $N^{j'}$}\label{sec:inaudibilidadNR}

In this Section we prove our main result.

\begin{theorem}\label{theo:inaudibilidadNR}
One cannot determine from the eigenvalues of the Laplace-Beltrami operator if a closed Riemannian manifold is naturally reductive.
\end{theorem}

To prove that the naturally reductive property of a closed Riemannian manifold is inaudible, we use Theorem \ref{theo:caracterizacionNR} and the isospectral pair $N^j$ and $N^{j'}$. To use this theorem on both manifolds, we first need to check if there exist $C,D\in\mathbb{R}$ with $D\neq 2C$ such that $J=C\cdot\textup{Id}_\mathfrak{v}$ and $B=D\cdot\textup{Id}_\mathfrak{z}$. Using \eqref{eq:j^2-base}, \eqref{def:matriz-j} and \eqref{def:corchete-j} we obtain in $N^j$ that
$$
\begin{aligned}
J=&(j_{Z_{\ii}})^2+(j_{Z_\j})^2+(j_{Z_\k})^2=-2\ \textup{Id}_{\mathfrak{v}},\\
B(Z^\mathfrak{z})=&\left[X_{\ii},j_{Z^\mathfrak{z}}X_{\ii}\right]+\left[X_{\j},j_{Z^\mathfrak{z}}X_{\j}\right]+\left[X_{\k},j_{Z^\mathfrak{z}}X_{\k}\right]\\
&+\left[Y_{\ii},j_{Z^\mathfrak{z}}Y_{\ii}\right]+\left[Y_{\j},j_{Z^\mathfrak{z}}Y_{\j}\right]+\left[Y_{\k},j_{Z^\mathfrak{z}}Y_{\k}\right]\\
=&4\ Z^\mathfrak{z}.
\end{aligned}
$$
Analogously, using \eqref{eq:j^2-base}, \eqref{def:matriz-j'} and \eqref{def:corchete-j'} we obtain in $N^{j'}$
$$
\begin{aligned}
J=&(j'_{Z_{\ii}})^2+(j'_{Z_\j})^2+(j'_{Z_\k})^2=-2\ \textup{Id}_{\mathfrak{v}},\\
B(Z^\mathfrak{z})=&\left[X_{\ii},j'_{Z^\mathfrak{z}}X_{\ii}\right]'+\left[X_{\j},j'_{Z^\mathfrak{z}}X_{\j}\right]'+\left[X_{\k},j'_{Z^\mathfrak{z}}X_{\k}\right]'\\
&+\left[Y_{\ii},j'_{Z^\mathfrak{z}}Y_{\ii}\right]'+\left[Y_{\j},j'_{Z^\mathfrak{z}}Y_{\j}\right]'+\left[Y_{\k},j'_{Z^\mathfrak{z}}Y_{\k}\right]'\\
=&4\ Z^\mathfrak{z}.
\end{aligned}
$$

\begin{proposition}
$N^j$ is a naturally reductive Riemannian manifold.
\end{proposition}
\begin{proof}
Let's identify $\mathfrak{z}=\mathbb{R}^3$ as the pure quaternion space $\mathbb{H}^*=\text{Span}\{{\ii},\j,\k\}$, and define the following map
    $$
    \begin{aligned}
    \{\cdot,\cdot\}:\mathfrak{z}\times\mathfrak{z}&\to\mathfrak{z}\\
    (X,Y)&\mapsto\{X,Y\}:=\text{Im}(X\cdot \bar{Y}),
    \end{aligned}
    $$
where $\bar{Y}$ denotes the conjugate of $Y\in\mathbb{H}^*$. It is well known that $(\mathfrak{z},\{\cdot,\cdot\})\cong(\mathbb{R}^3,\times)$, where $\times$ denotes the usual cross product in $\mathbb{R}^3$. Thus $(\mathfrak{z},\{\cdot,\cdot\})$ is a Lie algebra.

Now, let's take 
$
X^\mathfrak{z}=(x_7,x_8,x_9)\in\mathfrak{z},
Y^\mathfrak{z}=(y_7,y_8,y_9)\in\mathfrak{z}.
$
A straightforward computation gives us:
$$
j_{X^\mathfrak{z}}j_{Y^\mathfrak{z}}-j_{Y^\mathfrak{z}}j_{X^\mathfrak{z}}=\left(
\begin{array}{cc}
 A&0\\0&A
\end{array}
\right)=j_{\{X^\mathfrak{z},Y^\mathfrak{z}\}}
$$
where
$$
A=\begin{pmatrix}
0 & x_8 y_7-x_7 y_8 & x_9 y_7-x_7 y_9  \\
 x_7 y_8-x_8 y_7 & 0 & x_9 y_8-x_8 y_9  \\
 x_7 y_9-x_9 y_7 & x_8 y_9-x_9 y_8 & 0 
\end{pmatrix}.
$$
Moreover,
$$
\{X^\mathfrak{z},\cdot\}=\begin{pmatrix}
 0&-x_9&x_8\\
 x_9&0&-x_7\\
 -x_8&x_7&0
\end{pmatrix}\in\mathfrak{so(z)}.
$$
Thus, taking $\tilde{T}_{X^\mathfrak{z}}Y^\mathfrak{z}=\{X^\mathfrak{z},Y^{\mathfrak{z}}\}$, every condition in Theorem \ref{theo:caracterizacionNR} is satisfied, and
$$
T_XY=\frac{1}{2}j_{X^\mathfrak{z}}Y^\mathfrak{v}-\frac{1}{2}j_{Y^\mathfrak{z}}X^\mathfrak{v}+\frac{1}{2}[X^\mathfrak{v},Y^\mathfrak{v}]+\{X^\mathfrak{z},Y^\mathfrak{z}\}
$$
is a naturally reductive homogeneous structure in $N^j$ for every $X,Y\in\mathfrak{n}(j)$.
\end{proof}

\begin{proposition}
The Riemannian manifold $N^{j'}$ is not naturally reductive.
\end{proposition}
\begin{proof}
Let's consider nonzero $X^\mathfrak{z}=(x_7,x_8,x_9),Y^\mathfrak{z}=(y_7,y_8,y_9)\in\mathfrak{z}\subset\mathfrak{n}(j')$, and a map $\tilde{T}:\mathfrak{z}\times\mathfrak{z}\to\mathfrak{z}$ such that $\tilde{T}_{X^\mathfrak{z}}Y^\mathfrak{z}=(t_7,t_8,t_9)$. Then, by \eqref{def:matriz-j'},
$$
j'_{\tilde{T}_{X^\mathfrak{z}}Y^\mathfrak{z}}=
\begin{pmatrix}
 0 & 0 & 0 & 0 & -t_9 & t_8 \\
 0 & 0 & 0 & t_9 & 0 & -t_7 \\
 0 & 0 & 0 & -t_8 & t_7 & 0 \\
 0 & -t_9 & t_8 & 0 & 0 & 0 \\
 t_9 & 0 & -t_7 & 0 & 0 & 0 \\
 -t_8 & t_7 & 0 & 0 & 0 & 0
\end{pmatrix}.
$$
On the other hand,
$$
j'_{X^\mathfrak{z}}j'_{Y^\mathfrak{z}}-j'_{Y^\mathfrak{z}}j'_{X^\mathfrak{z}}=\left(
\begin{array}{cc}
 A&0\\0&A
\end{array}
\right)
$$
where
$$
A=\begin{pmatrix}
0 & x_8 y_7-x_7 y_8 & x_9 y_7-x_7 y_9  \\
 x_7 y_8-x_8 y_7 & 0 & x_9 y_8-x_8 y_9  \\
 x_7 y_9-x_9 y_7 & x_8 y_9-x_9 y_8 & 0 
\end{pmatrix}.
$$
Thus, none homogeneous structure of $N^{j'}$ could satisfy the relation \eqref{eq:condicion2-theorem} of Theorem \ref{theo:caracterizacionNR}. Therefore, none homogeneous structure of $N^{j'}$ could be naturally reductive.
\end{proof}

The previous propositions imply Theorem \ref{theo:inaudibilidadNR} due to $N^j$ and $N^{j'}$ are isospectral and non locally isometric Riemannian manifolds by Proposition \ref{prop:NjNj'isospectral} and \ref{prop:NjNj'noisometricas}. In other words, we have found an isospectral pair of closed Riemannian manifolds such that one of them is naturally reductive while its pair is not. This means that the eigenvalues of the Laplace-Beltrami operator cannot determine if a closed Riemannian manifold is naturally reductive.

\textbf{Author contributions:} All authors contributed equally to this research and in writing the paper.

\textbf{Funding:} The authors are supported by the grants GR21055 and IB18032 funded by Junta de Extremadura and Fondo Europeo de Desarrollo Regional. 
The first author is also partially supported by grant PID2019-10519GA-C22 funded by AEI/10.13039/501100011033 and by the grant GR24068 funded by Junta de Extremadura and Fondo Europeo de Desarrollo Regional.

\textbf{Conflicts of Interest:} The authors declare no conflict of interest. The founders had no role in the design of the study; in the collection, analyses, or interpretation of data; in the writing of the manuscript, or in the decision to publish the results.



\begin{thebibliography}{99}

\bibitem{AF.19}Arias-Marco, T. \& Fernández-Barroso, J. Inaudibility of k-D'Atri properties. {\em Symmetry}. \textbf{11}, 1316 (2019), \url{https://doi.org/10.3390/sym11101316}
\bibitem{AFF.15}Agricola, I., Ferreira, A. \& Friedrich, T. The classification of naturally reductive homogeneous spaces in dimensions $n\leq 6$. {\em Differential Geom. Appl.}. \textbf{39} pp. 59-92 (2015), \url{https://doi.org/10.1016/j.difgeo.2014.11.005}
\bibitem{AS.10}Arias-Marco, T. \& Schueth, D. On inaudible curvature properties of closed Riemannian manifolds. {\em Ann. Global Anal. Geom.}. \textbf{37}, 339-349 (2010), \url{https://doi.org/10.1007/s10455-009-9189-1}
\bibitem{AS.17}Arias-Marco, T. \& Schueth, D. Inaudibility of sixth order curvature invariants. {\em Rev. R. Acad. Cienc. Exactas F\'is. Nat. Ser. A Mat. RACSAM}. \textbf{111}, 547-574 (2017), \url{https://doi.org/10.1007/s13398-016-0311-5}
\bibitem{AS.58}Ambrose, W. \& Singer, I. On homogeneous Riemannian manifolds. {\em Duke Mathematical Journal}. \textbf{25}, 647 - 669 (1958), \url{https://doi.org/10.1215/S0012-7094-58-02560-2}
\bibitem{CFGV.23}Calviño-Louzao, E., Ferreiro-Subrido, M., Garc\'ia-R\'io, E. \& V\'azquez-Lorenzo, R. Homogeneous Riemannian structures in dimension three. {\em Revista De La Real Academia De Ciencias Exactas, F\'isicas Y Naturales. Serie A. Matem\'aticas}. \textbf{117}, 70 (2023)
\bibitem{CC.19}Calvaruso, G., López, M. \& Others Pseudo-Riemannian homogeneous structures. (Springer,2019)
\bibitem{E.94}Eberlein, P. Geometry of 2-step nilpotent groups with a left invariant metric. {\em Ann. Sci. École Norm. Sup. (4)}. \textbf{27}, 611-660 (1994).
\bibitem{GGSWW.98}Gordon, C., Gornet, R., Schueth, D., Webb, D. \& Wilson, E. Isospectral deformations of closed Riemannian manifolds with different scalar curvature. {\em Ann. Inst. Fourier (Grenoble)}. \textbf{48}, 593-607 (1998).
\bibitem{G.85}Gordon, C. Naturally reductive homogeneous Riemannian manifolds. {\em Canad. J. Math.}. \textbf{37}, 467-487 (1985), \url{https://doi.org/10.4153/CJM-1985-028-2}
\bibitem{G.78}Gray, A. Einstein-like manifolds which are not Einstein. {\em Geom. Dedicata}. \textbf{7}, 259-280 (1978), \url{https://doi.org/10.1007/BF00151525}
\bibitem{GW.97}Gordon, C. \& Wilson, E. Continuous families of isospectral Riemannian metrics which are not locally isometric. {\em J. Differential Geom.}. \textbf{47}, 504-529 (1997), \url{http://projecteuclid.org/euclid.jdg/1214460548}
\bibitem{K.66}Kac, M. Can one hear the shape of a drum?. {\em Amer. Math. Monthly}. \textbf{73}, 1-23 (1966), \url{https://doi.org/10.2307/2313748}
\bibitem{K.81}Kaplan, A. Riemannian nilmanifolds attached to Clifford modules. {\em Geom. Dedicata}. \textbf{11}, 127-136 (1981), \url{https://doi.org/10.1007/BF00147615}
\bibitem{KV.83}Kowalski, O. \& Vanhecke, L. Four-dimensional naturally reductive homogeneous spaces. {\em Rend. Sem. Mat. Univ. Politec. Torino}. pp. 223-232 (1984) (1983), Conference on differential geometry on homogeneous spaces (Turin, 1983)
\bibitem{KV.84}Kowalski, O. \& Vanhecke, L. A generalization of a theorem on naturally reductive homogeneous spaces. {\em Proc. Amer. Math. Soc.}. \textbf{91}, 433-435 (1984), \url{https://doi.org/10.2307/2045317}
\bibitem{KV.85}Kowalski, O. \& Vanhecke, L. Classification of five-dimensional naturally reductive spaces. {\em Math. Proc. Cambridge Philos. Soc.}. \textbf{97}, 445-463 (1985), \url{https://doi.org/10.1017/S0305004100063027}
\bibitem{L.98.2}Lauret, J. Naturally reductive homogeneous structures on 2-step nilpotent Lie groups. {\em Rev. Un. Mat. Argentina}. \textbf{41}, 15-23 (1998)
\bibitem{L.99}Lauret, J. Modified H-type groups and symmetric-like Riemannian spaces. {\em Differential Geom. Appl.}. \textbf{10}, 121-143 (1999), \url{https://doi.org/10.1016/S0926-2245(99)00002-9}
\bibitem{S.08}Schueth, D. Integrability of geodesic flows and isospectrality of Riemannian manifolds. {\em Math. Z.}. \textbf{260}, 595-613 (2008), \url{https://doi.org/10.1007/s00209-007-0290-5}
\bibitem{S.20}Storm, R. The classification of 7- and 8-dimensional naturally reductive spaces. {\em Canad. J. Math.}. \textbf{72}, 1246-1274 (2020), \url{https://doi.org/10.4153/s0008414x19000300}
\bibitem{Sz.99}Szabó, Z. Locally non-isometric yet super isospectral spaces. {\em Geom. Funct. Anal.}. \textbf{9}, 185-214 (1999), \url{https://doi.org/10.1007/s000390050084}
\bibitem{TV.83}Tricerri, F. \& Vanhecke, L. Homogeneous structures on Riemannian manifolds. (Cambridge University Press, Cambridge,1983), \url{https://doi.org/10.1017/CBO9781107325531}
\bibitem{W.82}Wilson, E. Isometry groups on homogeneous nilmanifolds. {\em Geom. Dedicata}. \textbf{12}, 337-346 (1982), \url{https://doi.org/10.1007/BF00147318}

\end{thebibliography}
\end{document}